%% file: cubegraph.tex
\documentclass[twoside,12pt]{article}
\usepackage{a4wide,amsmath,amssymb,epsfig}

\newcommand{\R}{{\mathbb R}}
\newcommand{\conv}{\mbox{\rm conv}}

\newcommand{\sse}{\subseteq}
\newcommand{\nn}{\mbox{\boldmath{$n$}}}
\newcommand{\xx}{\mbox{\boldmath{$x$}}}
\newcommand{\vv}{\mbox{\boldmath{$v$}}}

\newcommand{\sm}{{\setminus}}
\newcommand{\wt}{\widetilde}
\newtheorem {Theorem}{Theorem}
\newtheorem {Proposition}[Theorem]{Proposition}
\newtheorem {Conjecture}[Theorem]{Conjecture}
\newtheorem {Corollary}[Theorem]{Corollary}

\newtheorem {Lemma}[Theorem]{Lemma}

\def\proof{\par\medbreak\noindent{\bf Proof.}\enspace}
\def\qed{\vbox{\hrule
  \hbox{\vrule\hbox to 5pt{\vbox to 8pt{\vfil}\hfil}\vrule}\hrule}}
\def\endproof{\unskip \nobreak \hskip0pt plus 1fill \qquad \qed
\medskip\noindent}

\parskip=\smallskipamount

\title{\bfseries Neighborly cubical polytopes\thanks{%
Our work was supported by a DFG Gerhard-Hess-Forschungsf\"orderungspreis (Zi 475/1-2)
and by the German Israeli Foundation (G.I.F.) grant
I-0309-146.06/93.}}
\author{{\LARGE Michael Joswig} \quad { } \quad {\LARGE G\"unter M. Ziegler}\\[1mm]
Fachbereich Mathematik, MA~7-1\\ 
Technische Universit\"at Berlin\\
10623 Berlin, Germany\\
{\tt [joswig,ziegler]@math.tu-berlin.de}}
\date{March 23, 1999}

\begin{document}
\maketitle

\begin{abstract}
\noindent
Neighborly cubical polytopes exist: for any $n\ge d\ge 2r+2$, there is
a cubical convex $d$-polytope $C_d^n$ whose $r$-skeleton is
combinatorially equivalent to that of the $n$-dimensional cube.
This solves a problem of Babson, Billera \& Chan.

Kalai conjectured that the boundary~$\partial C_d^n$ of a neighborly
cubical polytope $C_d^n$ maximizes the $f$-vector among all cubical
$(d-1)$-spheres with $2^n$~vertices.  While we show that this is
true for polytopal spheres if $n\le d+1$, we also give a 
counter-example for $d=4$ and~$n=6$.

Further, the existence of neighborly cubical polytopes shows that the
graph of the $n$-dimensional cube, where $n\ge5$, is ``dimensionally
ambiguous'' in the sense of Gr\"unbaum.  We also show that the graph
of the $5$-cube is ``strongly $4$-ambiguous''.

In the special case $d=4$, neighborly cubical polytopes
have $f_3=\frac{f_0}4\log_2\frac{f_0}4$ vertices, so the
facet-vertex ratio $f_3/f_0$ is not bounded; this solves a problem
of Kalai, Perles and Stanley studied by Jockusch.
\end{abstract}

\section{Introduction.}

In Chapter~12 of his famous book~\cite{Gr} Gr\"unbaum discusses the
concept of {\em $k$-equivalence\/} of polytopes.  A $d$-polytope~$P$
is {\em $k$-equivalent\/} to a $d'$-polytope~$P'$ if the $k$-skeleta
of~$P$ and~$P'$ are combinatorially equivalent.  An interesting case
occurs when $\dim P\neq\dim P'$.  In this situation Gr\"unbaum calls
the $k$-skeleton~${\cal S}$ of either polytope {\em dimensionally
  ambiguous}.  Assume $d<d'$.  Then ${\cal S}$ is called {\em strongly
  $d$-ambiguous\/} if there is another $d$-polytope~$Q$, not
combinatorially equivalent to~$P$, whose $k$-skeleton is also
combinatorially equivalent to~${\cal S}$.  From the existence of
neighborly polytopes, such as the cyclic polytopes, it follows that
the $(\lfloor\frac{d}{2}\rfloor-1)$-skeleton of the $d$-simplex is
dimensionally ambiguous for $d\ge 5$.

Motivated by a problem of Kalai~\cite[Problem~3.4(iv)$^*$]{Z}, we
first constructed a $4$-polytope with the graph of the $5$-cube, thus
showing that the graph of the $5$-dimensional cube is dimensionally
ambiguous.  This turns out to be a special case of a $d$-polytope with
the $r$-skeleton of the $(d+1)$-cube, for
$r=\lfloor\frac{d}2\rfloor-1$.  We describe two constructions in
Section~\ref{sec:d+1}.  The resulting polytopes are cubical, and hence
they appear in the combinatorial classification of the cubical
polytopes with $2^{d+1}$~vertices due to Blind \&
Blind~\cite{BlBl0,BlBl2}.

More generally, in this paper we construct {\em neighborly cubical
  polytopes\/} in the sense of Babson, Billera \& Chan \cite{BBC}: for
every $n\ge d\ge2$ there is a cubical $d$-polytope with the
$r$-skeleton of the $n$-cube, for $r=\lfloor\frac d2\rfloor-1$.  In
particular, this yields a $4$-polytope with the graph of the $n$-cube
for every $n\ge4$.  The neighborly cubical polytopes~$C_d^n$ are
constructed as linear projections of ``deformed'' cubes, see
Section~\ref{sec:general}.  The combinatorics of the moment curve is
involved in an essential way; this is reminiscent of the construction
of neighborly (simplicial) polytopes.  We also give an explicit
combinatorial description of~$C_d^n$ which can be seen as a ``cubical
Gale evenness criterion.''

A result of Gr\"unbaum~\cite[12.2.1]{Gr} implies that no $d$-polytope
can have the $r$-skeleton of the $n$-cube for $r>\lfloor\frac
d2\rfloor-1$ for $n>d$; this is reviewed in Section~\ref{sec:ncp}.

Our construction specializes to a known phenomenon for $d=2$: There
are $n$-cubes whose $2$-dimensional ``shadows'' have $2^n$ vertices
(``projections that preserve the $0$-skeleton''). These were first
constructed by Murty \cite{Mu} and more explicitly by Goldfarb
\cite{Gol2}. They amount to linear programs for which the shadow
boundary pivot rule takes an exponential number of steps. In Amenta \&
Ziegler \cite{AZ} the Goldfarb cubes were interpreted as a special
case of a construction of ``deformed cubes,'' and indeed the
neighborly cubical polytopes constructed in this paper are projections
of deformed cubes $C^n(\varepsilon)$ as well.

An interesting new phenomenon occurs in the case $d=4$: Jockusch
\cite{Jock} had constructed examples of cubical $4$-polytopes for
which the facet/vertex ratio $f_3/f_0$ was higher than previously
expected, namely arbitrarily close to~$5/4$.  The neighborly
cubical polytopes show that indeed the ratio $f_3/f_0$ is not bounded
for cubical $4$-polytopes: these polytopes have $f_0=2^n$ vertices and
$f_3=(n-2)2^{n-2}=f_0\log(f_0/4)/4$ facets.

Kalai's cubical upper bound conjecture \cite[Conj.~4.2]{BBC} claimed
that among all cubical \mbox{$(d-1)$}-spheres with $2^n$~vertices, the
boundaries of cubical neighborly polytopes simultaneously maximize all
components of the $f$-vector. In Section~\ref{subsec:c46} we give a
counter-example, but we prove the claim in the special case $n=d+1$
for {\em polytopal\/} $(d-1)$-spheres.

A major part of the research on cubical polytopes is guided by the
comparison with simplicial polytopes.  In this paper, we extend this
analogy by providing cubical analogs to the cyclic (neighborly
simplicial) polytopes.  But we also offer a surprising instance where
the cubical case differs from the simplicial case: While
even-dimensional neighborly polytopes are always simplicial, it is not
true that a $4$-polytope with the graph of the $5$-cube must
necessarily be cubical --- we construct an explicit example in
Section~\ref{subsec:notcubical}.  Together with the existence of the
unique cubical $4$-polytope with the graph of the $5$-cube this
implies that the graph of the $5$-cube is strongly $4$-ambiguous.

\section{Neighborly cubical polytopes.}\label{sec:ncp}

We refer to \cite{Gr,Z} for general introductions to polytopes and
polytopal complexes. Two polytopes or polytopal complexes are {\em
  combinatorially equivalent\/} if their posets of faces are
isomorphic. In the following, a {\em $d$-cube\/} is any polytope that
is combinatorially isomorphic to the standard $d$-cube
$C_d=[0,1]^d\sse\R^d$.  A {\em combinatorial cube\/} is such a
$d$-cube, for any~$d$.  A {\em cubical polytope\/} is any polytope all
of whose proper faces are combinatorial cubes.  The 
{\em $k$-skeleton\/} of a polytope is the polytopal complex given by all
faces of dimension $k$ or less.  We say that $P$ 
{\em has the $k$-skeleton of a cube\/} 
if its $k$-skeleton is combinatorially equivalent to that
of a combinatorial cube.  A {\em neighborly cubical polytope\/} is a
cubical $d$-polytope (with $2^n$ vertices for some $n\ge d$) which has
the $(\lfloor\frac{d}2\rfloor-1)$-skeleton of a cube.  This notion was
introduced in \cite{BBC}, where neighborly cubical spheres were
constructed, and the question about the existence of neighborly
cubical polytopes was raised.  We start with explaining the choice of
parameters in this definition.

\begin{Proposition} {\bf (Characterization of Cubes \cite{BlBl0})}\\
Any cubical $d$-polytope has at least $2^d$ vertices.\\
If a cubical $d$-polytope
has exactly $2^d$ vertices, then it is a combinatorial $d$-cube.
\end{Proposition}

\begin{Corollary}~~\\
If all the $k$-faces of a $d$-polytope 
have $2^k$ vertices, for all $0\le k\le d-1$,
then the polytope is cubical.\\
If in addition the polytope has $2^d$ vertices, then it is
a combinatorial cube.
\end{Corollary}

It is well-known that the $f$-vector of a cubical polytope is subject
to restrictions that are similar to the Dehn-Sommerville equations
for simplicial/simple polytopes.

\begin{Proposition}\label{prop:dehn} {\bf (Cubical Dehn-Sommerville Equations
    \cite[9.4.1]{Gr})}\\
  Let $(f_0,\ldots,f_{d-1})$ be the $f$-vector of a cubical $d$-polytope.
  Then, for $0 \le k \le d-2$, 
\[
  \sum_{i=k}^{d-1}(-1)^i 2^{i-k} \binom{i}{k} f_i
  \ \ =\ \ (-1)^{d-1}f_k.
\]
\end{Proposition}

\begin{Lemma}\label{lem:simple} 
{\bf (Simple cubical polytopes \cite[Exercise~0.1, p.~23]{Z})}\\
Every  simple cubical $d$-polytope with $d>2$ is a $d$-cube.\\
(Every $2$-polytope is simple and cubical.)
\end{Lemma}

%

\begin{Proposition}\label{prop:double}~\\
  If a $d$-polytope $P$ has the $r$-skeleton of the $n$-dimensional
  cube, then all $k$-faces of~$P$ are cubes for $k\le 2r$.
\end{Proposition}

\proof Let $F$ be a $k$-face of~$P$ and $k\le 2r$.  By induction
on~$k$ we can assume that $F$ is cubical.  If $F$ is simple then $F$
is a cube by Lemma~\ref{lem:simple}.  Thus assume that $F$ has a
vertex~$v$ of degree~$k'>k$.  Let $a_1,\ldots,a_{k+1}$ be $k+1$
distinct vectors such that $v+a_i$ is a neighbor of~$v$ in~$F$.  As
$\dim F=k$, the vectors $a_1,\ldots,a_{k+1}$ are linearly dependent,
i.~e.\ we can choose $\lfloor\frac{k+1}2\rfloor$ vectors among them which
do not span a proper face of~$F$.  But 
$\lfloor\frac{k+1}2\rfloor \le \lfloor\frac{2r+1}2\rfloor = r$.  
This contradicts the assumption that $P$
has the $r$-skeleton of a cube.  \endproof

We note here that in the simplicial case more is true:
If a $d$-polytope $P$ has the $r$-skeleton of 
the $n$-dimensional {\em simplex\/},
then all $k$-faces of~$P$ are simplices for $k\le 2r+1$. In particular,
a $d$-polytope $P$ has the $r$-skeleton of the $n$-dimensional simplex
for $r\ge \lceil\frac{d}2\rceil$, then $P$ is a $d$-simplex.


\begin{Corollary}\label{cor:half}~\\
  If a $d$-polytope $P$ has the $r$-skeleton of the $n$-dimensional
  cube for $r\ge\frac{d}2$, then $P$ is a $d$-cube.
\end{Corollary}

But this is not good enough to establish the analogy to the simplicial
case. What if $P$ is a $2k$-polytope with the $(k-1)$-skeleton of the
$(2k+1)$-cube?  By Proposition~\ref{prop:double} all $(2k-2)$-faces
are cubes, but what about the $(2k-1)$-faces, i.~e.\ the facets?  It
turns out that the result of Proposition~\ref{prop:double} is sharp in
the sense that there are $2k$-polytopes which have the
$(k-1)$-skeleton of a $(2k+1)$-cube but which are {\em not\/} cubical:
In Section~\ref{subsec:notcubical} we present an example of a non-cubical
$4$-polytope with the graph of the $5$-cube.

The proof of the following is based on a theorem of van Kampen and
Flores, see~\cite[11.1.3 and 11.3]{Gr}.

\begin{Proposition}%
  {\bf (Gr\"unbaum~\cite[11.2.1]{Gr})}\\
  Let $P$ be a $d$-polytope with the $r$-skeleton of an $n$-polytope, with
  $n>d$.  Then $r\le \lfloor \frac{d}2\rfloor-1$.
\end{Proposition}

\section{Projections.}\label{sec:projections}

We discuss the effect of orthogonal projections on polytopes.
Everything in this section is well-known; it is included for the sake
of completeness.

Let $P$ be a full-dimensional polytope.  A vector~$\nn$ is  {\em
  normal\/} with respect to a facet~$F$ of~$P$ if it is orthogonal
to~$F$ and it points ``to the outside,'' that is, if the linear functional
corresponding to~$\nn$ and restricted to~$P$ attains its maximum at
the points in~$F$.  A vector is called {\em normal\/} with respect to
a face~$G$ if it is a positive linear combination of normal vectors of
all facets of~$P$ containing~$G$.  Equivalently, the linear
functional corresponding to a normal vector of~$G$ attains its maximum
at the points in~$G$.  Obviously, every facet has a unique normal vector
of length~$1$, while a face of higher codimension does not.

Consider an orthogonal projection~$\pi$ onto some proper affine subspace.

\begin{Lemma}\label{lem:orthogonal}~\\
  If the face~$G$ has a normal vector orthogonal to the direction
  of projection, then the image $\pi(G)$ is a face of the
  polytope~$\pi(P)$.  Conversely, if $\overline G$ is a face of~$\pi(P)$,
  then the full preimage~$\pi^{-1}(\overline{G})$ is a face of~$P$ with a
  normal vector orthogonal to the direction of projection.
\end{Lemma}

\begin{figure}[htbp]
  \begin{center}
    \epsfig{file=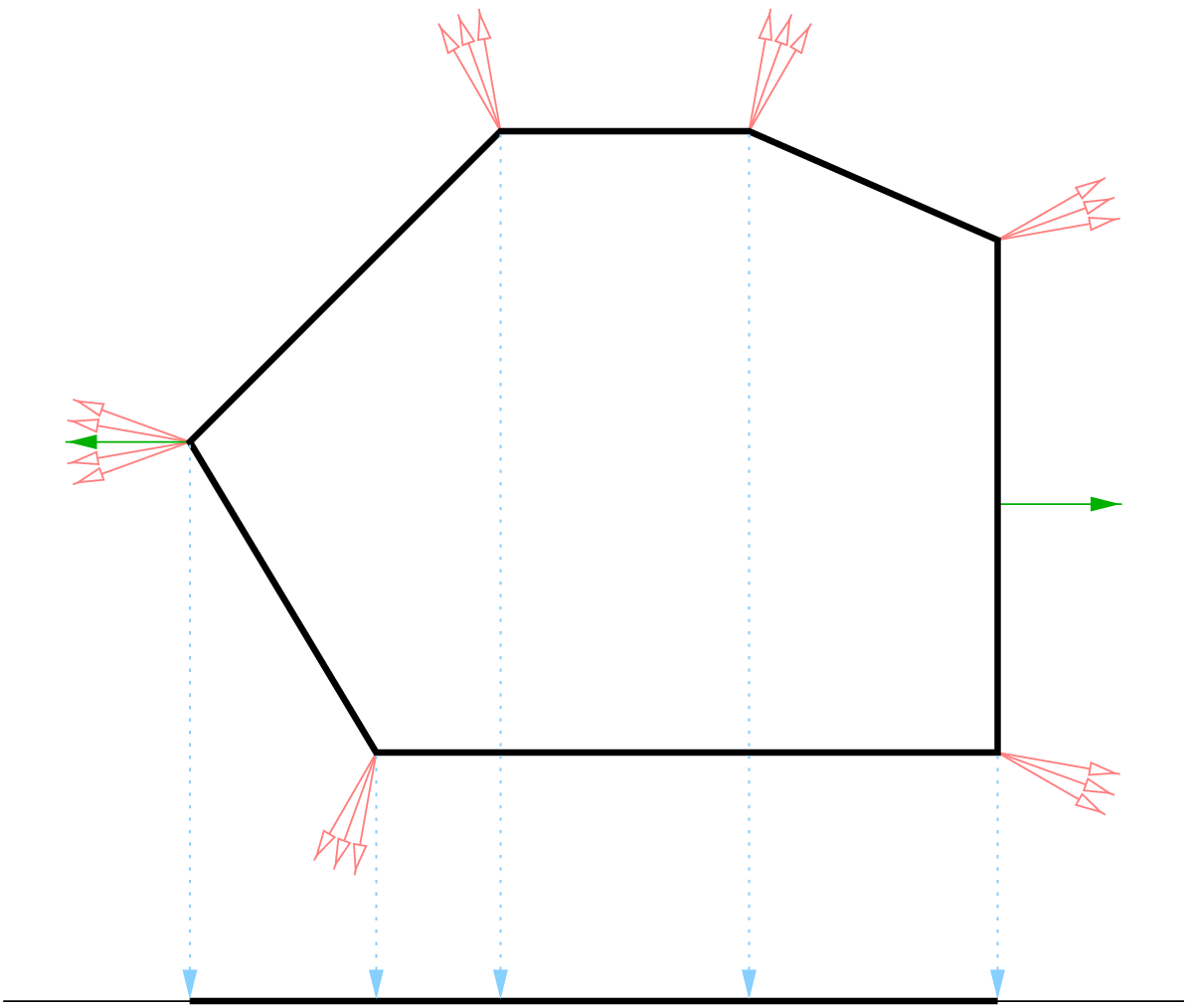,width=6cm}
    \caption{Orthogonal projection of a $2$-polytope.}
    \label{fig:projection}
  \end{center}
\end{figure}

The lemma above characterizes the ``shadow boundary,''
that is, it describes the faces that are mapped to the boundary
of the projection.
Note that this may include faces of~$G$ that are mapped to faces of
lower dimension, cf.\
Figure~\ref{fig:projection}.

\begin{Lemma}~~\\
  The restriction of $\pi$ to a face~$G$ is injective if and only if
  $G$ has a normal vector which not orthogonal to the direction of
  projection.
\end{Lemma}

Combining the two lemmas above gives the following characterization.

\begin{Corollary}\label{cor:normal}~\\
  A face~$G$ is mapped onto a face~$\pi(G)$ of the same dimension
  if and only if it has a normal vector which is orthogonal to the 
  direction of projection and another one which is not.
\end{Corollary}

We can also say something about the shape of the projection in this case.

\begin{Lemma}\label{lem:isomorphic}~\\
  A face~$G$ is mapped onto a face~$\pi(G)$ of the same dimension if
  and only if the faces $G$ and~$\pi(G)$ are affinely isomorphic.
\end{Lemma}

\section{The case \ \boldmath{$n=d+1$}.}\label{sec:d+1}

\subsection{First construction.}

Let $d=2r$ be even, $Q:=[-1,+1]^r$, and 
\[
\wt{P}\ \ :=\ \ \conv(Q\times 2Q\times\{-1\}\ \cup\ 2Q\times
Q\times\{1\})\ \ \sse\ \ \R^{d+1}.
\]
This is clearly a combinatorial $(d+1)$-cube,
with the complete the linear description 
\begin{eqnarray*}
\wt P\ \ =\ &\{\tbinom{\xx}{x_{d+1}}\in\R^{d+1}: &
-1\le x_{d+1}\le 1, \\
&&\pm\,2x_i\ \le\ 3-x_{d+1}\textrm{\ \quad for } 1\le i\le r,\\
&&\pm\,2x_i\ \le\ 3+x_{d+1}\textrm{\ \quad for } r<   i\le d \ \ \}.
\end{eqnarray*}
The projection $\pi:\R^{d+1}\longrightarrow\R^d$
that deletes the last coordinate yields the $d$-polytope
\[
P\ \ :=\ \ \pi(\wt{P})\ \ =\ \ \conv(Q{\times}2Q\ \cup\ 2Q{\times}Q)\sse\R^d.
\]
One Fourier-Motzkin elimination step \cite[Sect.~1.2]{Z}
shows that $P$ can also be described in terms of its facets by 
\begin{eqnarray*}
P&=&\{\xx\in\R^d: \qquad\ \pm\, x_i \le 2 
\textrm{\quad for } 1\le i\le 2r,\\
 & &\qquad\qquad\ \,\pm\, x_i \pm x_j \le 3 
\textrm{\quad for } 1\le i\le r<j\le 2r    \ \ \}.
\end{eqnarray*}
This $P$ is a cubical $2r$-polytope with
$f_0=2^{d+1}$ vertices and $f_{d-1}=4r+(2r)^2=d(d+2)$ facets.

It is also easy to see that $P$ has the $(r-1)$-skeleton of a
$(d+1)$-cube, either using the criteria of Section~\ref{sec:projections},
or by direct verification from the complete description 
of the polytope in terms of facets {\em and\/} vertices.

This example of a polytope with the $\frac{d}2$-skeleton of a cube
is amazing because of its simplicity, and also because of its
symmetry: It has a vertex-transitive symmetry group,
and only two orbits of facets.

\subsection{Second construction.}

Blind \& Blind completed a classification of the (combinatorial types)
of cubical $d$-polytopes with $2^{d+1}$ vertices in~\cite{BlBl2}.
{}From this classification, we derive below that for even~$d$, there
is exactly one combinatorial type of a $d$-polytope with the
$\frac{d}2$-skeleton of the $(d+1)$-cube; for odd $d$, there are
precisely two combinatorial types.  The description of the polytopes
given by Blind \& Blind also implies the following construction given
above.

A (cubical) $d$-polytope~$P$ whose boundary complex~$\partial P$ is isomorphic
to a subcomplex of the $d$-skeleton of some (higher-dimensional) cube
is called {\em liftable}.  We get all distinct
combinatorial types of liftable $d$-polytopes (and
$(d-1)$-spheres) with at most $2^{d+1}$ vertices as follows:

\begin{quote}
  The cube $C^{d+1}$ has $d+1$ pairs of opposite facets $F^{\pm}_i$
  $(1\leq i\leq d+1)$.  For $k,m\ge 0$ and $l\ge 1$ with $k+l+m=d+1$
  let $B(k,l,m)$ be the cubical $d$-ball in the boundary of $C^{d+1}$
  formed by $F^{\pm}_i$ $(1\leq i\leq k)$ and $F^+_i(k+1\leq i\leq
  k+l)$.
\end{quote}

The combinatorial types $P(k,l,m)$ of liftable $d$-polytopes with at
most~$2^{d+1}$ vertices are given by the boundary complexes of the
cubical balls $B(k,l,m)$.  The number of vertices of~$P(k,l,m)$
equals~$2^{d+1}$ if and only if $k,m\ge 1$.  Note that $P(k,l,m)$ is
combinatorially equivalent to~$P(m,l,k)$; thus in the following let
$k\ge m$.  There are $\lfloor d^2/4\rfloor$ suitable triples $(k,l,m)$
with $k\ge m\ge 1$.

\begin{Theorem}\label{Thm:BBclass}
  {\bf (Classification for \boldmath{$n=d+1$}, Blind \& Blind \cite[Theorem~3]{BlBl2})}\\
  For $d\ge4$, all the combinatorial types of cubical $d$-polytopes with
  $2^{d+1}$ vertices are given by the liftable polytopes $P(k,l,m)$
  with $k,l,m\ge1$, $k+l+m=d+1$ and $k\ge m$, plus in addition the
  ``$2$-fold non-linearly capped'' $d$-polytope $P^d_{NLC}$.
\end{Theorem}

A polytope $P(k,l,m)$ is {\em $r$-neighborly\/} if for every
$r$-face~$H$ of $C^{d+1}$, there is some facet of~$B(k,l,m)$
containing~$H$, and some facet of the complement~$C^{d+1}\setminus
B(k,l,m)$ containing~$H$ as well.  Every $r$-face of $C^{d+1}$ is
contained in exactly $d+1-r$ facets of $C^{d+1}$. Thus $P(k,l,m)$ is
$r$-neighborly if and only if $k+l\leq (d+1-r)-1$ and $l+m\leq
(d+1-r)-1$, that is $r+1\le m$ and $r+1\le k$, where the first
condition implies the second because of $k\geq m$.

Now we consider the case $r=\lfloor \frac{d}{2}\rfloor -1$.

\begin{Corollary}~~\\
  If $d$ is even, then $P(\frac{d}{2},1,\frac{d}{2})$ is the unique
  $(\frac{d}{2}-1)$-neighborly cubical polytope with
  $2^{d+1}$~vertices.  If $d$ is odd, then there are precisely two
  $(\frac{d-1}{2}-1)$-neighborly cubical polytopes with
  $2^{d+1}$~vertices, namely $P(\frac{d+1}{2},1,\frac{d-1}{2})$ and
    $P(\frac{d-1}{2},1,\frac{d-1}{2})$.
\end{Corollary}

Now we compute the $f$-vector of $P(k,l,m)$, denoted by
$f(P(k,l,m))=(f_0,\ldots,f_{d-1})$. With the same reasoning as above,
an $i$-face of $C^{d+1}$ is a face of $P(k,l,m)$ if it lies in some
facet of $B(k,l,m)$, and also in some facet not in $B(k,l,m)$. We
deduce that
\[
f_{d+1-i}\ \ =\ \ \binom{d+1}{i}2^i\ -\ 
\sum_{j=0}^i\binom{l}{i-j}\Big\{\binom{k}{j}+\binom{m}{j}\Big\}2^j.
\]


\subsection{The Cubical Upper Bound Conjecture.}

From the analogy to the simplicial case one is tempted to expect that
the neighborly cubical polytopes achieve equality for the cubical
upper bound conjecture.

\begin{Conjecture}\label{Conj:CUBC} 
  {\bf (Cubical Upper Bound Conjecture, Kalai \cite[Conjecture 4.2]{BBC})}\\
   Let $P$ be a cubical $(d-1)$-sphere with $f_0(P)=2^n$ vertices. Then
   its number of facets is bounded by that of $C_d^n$, that is,
   $f_{d-1}(P)\le f(n,d)$.  Moreover,
\[
f_i(P)\leq f_i(C_d^n)\qquad\textrm{ for }1\leq i\leq d-1.
\]
\end{Conjecture}

\begin{Theorem}~\label{thm:n=d+1}~\\
  In the special case $n=d+1$, Conjecture~\ref{Conj:CUBC} is 
{\em true\/} if restricted to cubical polytopes.  
However, it is {\em
    false\/} for spheres even for $d=4$ and~$n=6$.
\end{Theorem}

\proof\\
{\bf (1)} 
The proof for polytopes relies on the classification of
Theorem~\ref{Thm:BBclass}. Here we can disregard the ``$2$-fold
non-linearly capped'' $d$-polytope $P^d_{NLC}$, since it has the same
$f$-vector as the ``$2$-fold linearly capped'' $d$-polytope, which is
$P(d-1,1,1)$, see \cite[Figure~1]{BlBl2}.
Thus our problem is to minimize, for fixed $d$ and~$i$, the function 
\[
\delta_i(k,l,m)\ \ =\ \ 
\sum_{j=0}^i\binom{l}{i-j}\Big\{\binom{k}{j}+\binom{m}{j}\Big\}2^j
\]
subject to the restrictions $k,l,m\ge1$, $k+l+m=d+1$ and $k\ge m$.
For this, we note the simple properties and inequalities
\begin{itemize}
\item[] $\delta_i(k,l,m)\ \ge\ \delta_i(k-1,l,m+1)$ \ \ for $k>m$,\\[.5mm]
        $\delta_i(k,l,m)\ \ge\ \delta_i(k+1,l-2,m+1)$ \ \ for $l\ge2$,\\[.5mm]
        $\delta_i(k,2,m)\ \ge\ \delta_i(k,1,m+1)$\ \ for $k>m$, and\\[.5mm]
        $\delta_i(k,2,k)\  = \ \delta_i(k+1,1,k)$,
\end{itemize}
from which the claim immediately follows.

{\bf (2)} The second claim is verified via an explicit construction
in Section~\ref{subsec:c46}.
\endproof
  
For even~$d$, the first part of this
result also follows from Babson, Billera \& Chan
\cite[Thm.~4.3]{BBC}, who used Adin's ``cubical $h$-vector''
\cite{Ad}.

\subsection{A non-cubical polytope.}\label{subsec:notcubical}

Here we give an example of a non-cubical $4$-polytope with the same
graph as the $5$-cube.  For this it is helpful to have yet another
coordinate representation of the cubical $4$-polytope with the
graph of the $5$-cube. 

Let $P$ be the polytope defined as the convex hull of the following
$32$~points in~$\R^4$.  Note that each row corresponds to four
distinct points due to arbitrary variation of the signs.
\[
\begin{array}{crrrrc}
  (&\pm 1 & \pm 1 & 1 & 1 &)\\
  (&\pm 1 & \pm 1 & 4 & 1 &)\\
  (&\pm 2 & \pm 2 & 3 & 4/5 &)\\
  (&\pm 2 & \pm 2 & 2 & 4/5 &)\\
  (&\pm 3 & \pm 3 & 2 & 1/2 &)\\
  (&\pm 3 & \pm 3 & 3 & 1/2 &)\\
  (&\pm 4 & \pm 4 & 0 & 0 &)\\
  (&\pm 4 & \pm 4 & 5 & 0 &)
 \end{array}
\]
All these points are vertices.  Moreover, the last eight vertices span
a facet~$F=\{\xx\in P:x_4=0\}$, which is a $3$-cube.  The graph of~$P$
is isomorphic to the graph of the $5$-cube.

\begin{figure}[htbp]
  \begin{center}
    \begin{minipage}[t]{7.5cm}
      \begin{center}
        \epsfig{file=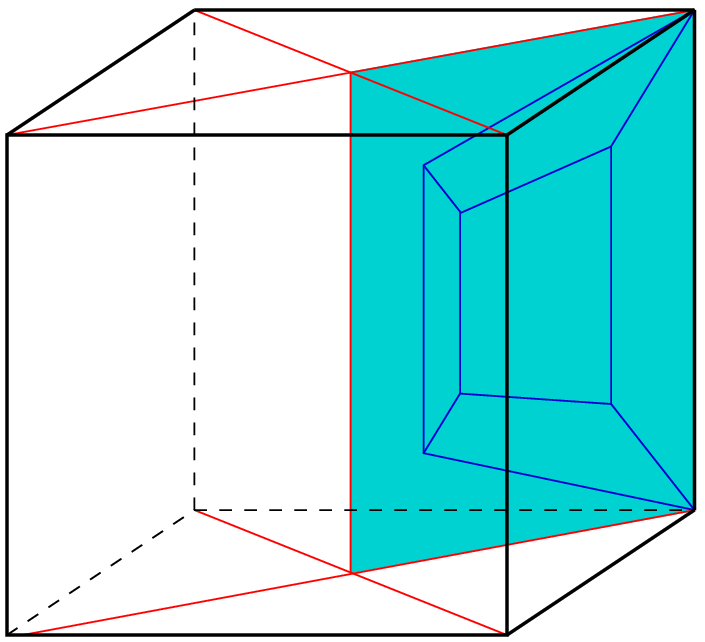,width=6cm}
      \end{center}
    \end{minipage}
    \begin{minipage}[t]{7.5cm}
      \begin{center} 
        \input{quarter.pstex_t}
      \end{center}
    \end{minipage}
    \caption[Schlegel diagram.]{Schlegel diagram and a section.}
    \label{fig:schlegel}
  \end{center}
\end{figure}

Projecting the polytope~$P$ onto the facet~$F$ from a point beyond~$F$
yields a polytopal complex, the {\em Schlegel diagram}\/ of~$P$ with
respect to~$F$, which is essentially equivalent to the boundary
complex of~$P$ \cite[Sect.~5.3]{Gr} \cite[\S5.2]{Z}.  For our example,
the whole Schlegel diagram has the same symmetry group of order~$16$
as a prism over a square.  Thus it is sufficient to consider a
diagonal section as indicated in Figure~\ref{fig:schlegel}.

\vskip-4mm

Now remove the $2$-face spanned by the $4$~vertices $(\pm 1,\pm
1,4,1)$ and merge the two facets containing it.  This yields a regular
cell complex~$C$ whose $1$-skeleton is the same as before, that is, it
is isomorphic to the graph of the $5$-cube.  Clearly, $|\!|C|\!|$ is
homeomorphic to the $3$-sphere.  But, as realized, $C$ is not
polytopal because the angle~$\alpha$ in Figure~\ref{fig:schlegel}
exceeds~$\pi$.

It is fairly obvious that the modified cell complex~$C$ can be
realized as a polytopal complex by an appropriate change of the
coordinates, such that finally $\alpha$ becomes less than~$\pi$.  It
may be less obvious that --- for a special choice of coordinates ---
the transformed diagram with $\alpha<\pi$ can be lifted to a
$4$-polytope~$P'$.  Here is a realization:
$$
\begin{array}{crrrrc}
  (&\pm 1 & \pm 1 & 1 & 0 &)\\
  (&\pm 2 & \pm 2 & 4 & 0 &)\\
  (&\pm 3 & \pm 3 & 3 & 1 &)\\
  (&\pm 3 & \pm 3 & 2 & 5/4 &)\\
  (&\pm 4 & \pm 4 & 2 & 21/20 &)\\
  (&\pm 4 & \pm 4 & 3 & 9/5 &)\\
  (&\pm 56/13 & \pm 56/13 & 0 & 779/260 &)\\
  (&\pm 5 & \pm 5 & 16 & 0 &)
 \end{array}
$$
The polytope~$P'$ has one facet that has $12$ vertices, namely $$(\pm
1,\pm 1,1,0),\ (\pm 2,\pm 2,4,0),\ (\pm 5,\pm 5,16,0).$$
Therefore
$P'$ is a $4$-polytope with the graph of the $5$-cube which is not
cubical.

\section{Alternating Oriented Matroids and Cyclic
  Polytopes.}\label{sec:om}

In Section~\ref{sec:general} we will construct a class of polytopes
which are ``cubical relatives'' of the cyclic polytopes.  The
combinatorial structure of the cyclic polytopes is well-known
\cite[Sect.~4.6]{Gr} \cite[Example~0.6]{Z}.  We give a brief account
in the framework of oriented matroids, which also captures the
``interior combinatorial structure'' of the cyclic polytopes.

Let $C_d(n)$ be the cyclic $d$-polytope on $n$~vertices.  It can be
realized as the convex hull of $n$~points on any curve of order~$d$,
such as the moment curve $t\mapsto(t,t^2,t^3,\ldots,t^d)$ in~$\R^d$.

Any point configuration in~$\R^d$, and thus any polytope via its
vertices, gives rise to an oriented matroid \cite[Sect.~6.4]{Z}
\cite{BLSWZ}.  The positive cocircuits of the oriented matroid
bijectively correspond to the facets of the polytope.

For a cyclic polytope~$C_d(n)$ the associated oriented matroid is
known as the {\em alternating oriented matroid\/}~${\cal C}(n,d+1)$ of
rank~$d+1$ on $n$~points \cite[Sects.~3.4 and~9.4]{BLSWZ}.
Homogenizing the vertices of~$C_d(n)$ it can be represented by the
(rows of the) matrix
\[
\left(
\begin{array}{ccccc}
1 & t_1 & t_1^2 & \cdots & t_1^d \\
1 & t_2 & t_2^2 & \cdots & t_2^d \\
\vdots & \vdots & \vdots & \ddots & \vdots \\
1 & t_n & t_n^2 & \cdots & t_n^d
\end{array}
\right)
\]
for $t_1<t_2<\cdots<t_n$.

There is a notion of duality for oriented matroids which generalizes
duality of projective spaces.  In the case of the alternating oriented 
matroids, the dual of~${\cal C}(n,d+1)$ is obtained
from~${\cal C}(n,n-d-1)$ by reorienting every other row in the
representation above, 
\[
{\cal C}(n,d+1)^*\ \ =\ \ {}_{\overline{\{2,4,6,\ldots\}}}{\cal C}(n,n-d+1),
\]
(see \cite[pp.~108-109]{BLV}, \cite[Sect.~2]{Stu}),
that is, by the rows of the $n\times (n-d+1)$-matrix
$$\left(
\begin{array}{rrrcr}
1 \hphantom{)^{n+1}}& t_1 & t_1^2 & ~~\cdots~~ & t_1^{n-d-1} \\
-1\hphantom{)^{n+1}} & -t_2 & -t_2^2 & \cdots & -t_2^{n-d-1} \\
1 \hphantom{)^{n+1}}& t_3 & t_3^2 & \cdots & t_3^{n-d-1} \\
-1\hphantom{)^{n+1}} & -t_4 & -t_4^2 & \cdots & -t_4^{n-d-1} \\
\vdots~ \hphantom{)^{n+1}}& \vdots~ & \vdots~ & \ddots & \vdots\qquad\\
{(-1)}^{n+1} & {(-1)}^{n+1}t_n & {(-1)}^{n+1}t_n^2 & \cdots & {(-1)}^{n+1}t_n^{n-d-1}
\end{array}
\right)$$

{}From the representation of the alternating matroid and its dual
given above it is obvious that that a deletion (omitting a row in the
primal oriented matroid) or a contraction of the first element
(omitting the first row and the first column in the primal oriented
matroid, which amounts to omitting the first row in the dual) again gives an
alternating matroid (on fewer points and, in the second case, of
smaller rank).

The number of facets of~$C_d(n)$, that is, the number of 
positive cocircuits of
the alternating matroid~${\cal C}(n,d+1)$, is known to be
$$\binom{n-\lceil\textstyle{\frac d2}\rceil}{\lfloor\textstyle{\frac
    d2}\rfloor}+\binom{n-1-\lceil \textstyle{\frac
    {d-1}{2}}\rceil}{\lfloor\textstyle{\frac {d-1}{2}}\rfloor},$$
by Gale's evenness criterion \cite[p.~63]{Gr} \cite[Thm.~0.7]{Z}.

\section{The general case.}\label{sec:general}

The following is our main theorem: ``Neighborly cubical polytopes
exist!''  --- they can be obtained as projections of deformed cubes.
As mentioned in the introduction, the case $d=2$ and the case $n=d+1$
were known previously.

\begin{Theorem}~~\\
  For any $n\geq d\geq 2r+2$, there exists a combinatorial $n$-cube
  $C^n\subseteq \R^n$ and a linear projection map $\pi:\R^n\rightarrow
  \R^d$ such that $C_d^n:=\pi(C^n)$ is a cubical $d$-polytope whose
  $r$-skeleton is isomorphic to that of $C^n$ (via $\pi$).
\end{Theorem}

\proof We first construct a combinatorial $n$-cube
$C^n(\varepsilon)\subseteq \R^n$ that depends on a parameter
$\varepsilon >0$; then we verify that for $\varepsilon$ sufficiently
small the projection $\pi:\R^n\rightarrow\R^d$ to the last $d$
coordinates preserves the $r$-skeleton; and finally we argue that
$\pi(C^n(\varepsilon))=C_d^n$ is cubical.

{\bf (A) } For $0<\varepsilon \leq 1$ define
$C^n(\varepsilon)\subseteq \R^n$ as the solution set of
\begin{equation}
\varepsilon|x_k|\ \ \ \leq\ \ \ \frac{2^{\binom{k}{2}}}
{\varepsilon^{k-1}}-(-1)^k\sum\limits^{k-1}_{j=1}
{\binom{k-2}{j-1}}x_j\qquad\textrm{ for }1\leq k\leq n.
\end{equation}

This set is a combinatorial $n$-cube. To see this, we verify by
induction on $k$ that all solutions to the first $k$ conditions
of (1) satisfy

\begin{equation}
|x_k|\ \ <\ \ \frac{2^{\binom{k}{2}+1}}{\varepsilon^k}.
\end{equation} 

In fact, the upper bound of (2) increases with $k$, and for $k=1$
we have $\varepsilon|x_1|\leq 1$, so
$|x_1|<\frac{2}{\varepsilon}$ is surely satisfied. Thus we may
use induction, and for $k\ge2$ estimate

$$
\begin{array}{lcl}
\varepsilon|x_k| & \leq &
\displaystyle{\frac{2^{\binom{k}{2}}}{\varepsilon^{k-1}}-(-1)^k
\sum\limits^{k-1}_{j=1} {\binom{k-2}{j-1}}x_j}\\
& < & 
\displaystyle{\frac{2^{\binom{k}{2}}}{\varepsilon^{k-1}}+
\sum\limits^{k-1}_{j=1} {\binom{k-2}{j-1}}
\frac{2^{\binom{k-1}{2}+1}}{\varepsilon^{k-1}}}\\
& \le & 
\displaystyle{\frac{2^{\binom{k}{2}}}{\varepsilon^{k-1}}+2^{k-2}
\frac{2^{\binom{k-1}{2}+1}}{\varepsilon^{k-1}}\ \ =\ \ 
\frac{2^{\binom{k}{2}+1}}{\varepsilon^{k-1}}}.
\end{array}$$

In this computation the second term is always smaller in absolute
value than the first, that~is,
$$\frac{2^{\binom{k}{2}}}{\varepsilon^{k-1}}-\left|
\sum\limits^{k-1}_{j=1} {\binom{k-2}{j-1}}x_j\right| \ \ >\ \  0\qquad
\textrm{for } 1\leq k\leq n,$$
and from this we see that $C^n(\varepsilon)$ is a combinatorial
cube. (It is an iterated deformed product in the sense of~\cite{AZ}).
\smallskip

{\bf (B) }  Now let $\pi:\R^n\rightarrow \R^d$ be the
projection to the last $d$ coordinates, $(x_1,\ldots,
x_n)\longmapsto (x_{n-d+1},\ldots, x_n).$

{\em{\bf Claim 1.} For sufficiently small $\varepsilon > 0$ and $n\geq
  d\geq 2r+2$, the orthogonal projection
  $\pi:C^n(\varepsilon)\longrightarrow
  \pi(C^n(\varepsilon))=:C_d^n(\varepsilon)$ preserves the
  $r$-skeleton. That is, $\pi$ restricts to an isomorphism
$$\pi:C^n(\varepsilon)^{[r]}\longrightarrow
C_d^n(\varepsilon)^{[r]}$$
of polytopal complexes.}

To verify this, we need to see that every $r$-face $F$ of
$C^n(\varepsilon)$ is mapped bijectively to an $r$-face $\pi(F)$
of $C^n(\varepsilon)$. By Corollary~\ref{cor:normal} this is
equivalent to the condition that for every $r$-face $F$ there is a normal
vector which is orthogonal to the direction of projection and another
one which is not.

In our specific situation, let $F$ be an $r$-face of
$C^n(\varepsilon)$, and let $\vv\in C^n(\varepsilon)$ be a
vertex of $F$. Then there are unique signs
$\sigma_1,\ldots,\sigma_n\in
\{+1,-1\}$ such that $\vv$ is determined by
$$
(-1)^kÊ\sum\limits^{k-1}_{j=1}{\binom{k-2}{j-1}}v_j +\sigma_k
\varepsilon v_k=
\frac{2^{\binom{k}{2}}}{\varepsilon^{k-1}}\qquad\textrm{ for
}1\leq k\leq n,
$$
while $F$ is characterized by the additional choice of a set
$S\in{\binom{[n]}{n-r}}$ of $n-r$ indices:
$$
F=\{{\xx}\in C^n(\varepsilon):(-1)^k
\sum\limits^{k-1}_{j=1}{\binom{k-2}{j-1}} x_j + \sigma_k
\varepsilon x_k = 
\frac{2^{\binom{k}{2}}}{\varepsilon^{k-1}}\ \ \textrm{ for all }
k\in S\},
$$
where $C^n(\varepsilon)$ itself is given by
$$
C^n(\varepsilon)=\{{\xx}\in\R^n:(-1)^k
\sum\limits^{k-1}_{j=1}{\binom{k-2}{j-1}} x_j + \varepsilon |x_k|
\leq \frac{2^{\binom{k}{2}}}{\varepsilon^{k-1}}\ \ \textrm{ for }
1\leq k\leq n\}.
$$
In order to show that all $\pi({\vv})$ are vertices of
$\pi(C^n(\varepsilon))$, we must thus check that, for any choice
$\Sigma=(\sigma_1,\ldots,\sigma_n)\in\{\pm 1\}^n$ of signs, the rows
of the $n\times(n-d)$-matrix
$$
A(\Sigma)=\left(
\begin{array}{rrrcr}
\pm\varepsilon&0              &              & &\\
 1             &\pm\varepsilon&0             & &\\
 -1            &-1            &\pm\varepsilon&~~\ddots~~&\\
 1             &2             &1             &\ddots&0\\
 -1            &-3            &-3            &\ddots     
&\pm\varepsilon\\
 1             &4             &&&\\
\vdots &&\ddots&&
\end{array}
\right)
$$
given by
$$
a_{kj}=\left\{
\begin{array}{lll}
0                       & \textrm{ for }& j>k\\
\sigma_k\varepsilon     & \textrm{ for }& j=k\\
(-1)^k{\binom{k-2}{j-1}}& \textrm{ for }& j<k
\end{array}
\right.
$$
have a positive linear dependence.

We also have to show that each $r$-face has a normal vector which is
not orthogonal to the direction of projection.  But, this is obvious:
for any set of rows of the matrix~$A(\Sigma)$ (for an arbitrary
vector~$\Sigma$ of signs) set, e.~g., the coefficient of the last row
to~$1$ and all the others positive but sufficiently small.  This
positive linear combination yields a non-zero vector.

At an $r$-face $F\subseteq C^n(\varepsilon)$ only $n-r$ restrictions
are tight, so Claim~1 now reduces to the following.

{\bf Claim 2.} {\em For sufficiently small $\varepsilon >0$ and for
  every choice $\Sigma\in\{\pm 1\}^n$ of signs, every set of $n-r$
  rows of $A(\Sigma)$ has a positive dependence.}

Let $\overline A(\Sigma):= (a_{kj})_{2\leq k\leq n, 1\leq j\leq
  n-d}\in \R^{(n-1)\times(n-d)}$ be obtained by deleting the first row
of~$A(\Sigma)$. An index set $S\subseteq\{2,3,\ldots,n\}$ will be
called {\em alternating\/} if it alternates between odd and even
numbers; for example $\{2,3,6,9\}$ and $\{3,4,5,8\}$ are alternating,
but $\{2,3,5,6\}$ is not. A set of rows of $S$ is {\em alternating\/}
if the corresponding index set is alternating. Using this concept, we
formulate the following Claim 3, which clearly implies Claim 2.

{\bf Claim 3.}\\
{\bf (i)} {\em If $\varepsilon=0$, then all maximal minors of
  $\overline A(\Sigma)$ have non-zero determinant. \\
  {\bf (ii)} If $\varepsilon >0$ is so small that all maximal minors
  of $\overline A(\Sigma)$ have the same sign as for $\varepsilon=0$,
  then every alternating subset of $n-d+1$~rows of~$\overline
  A(\Sigma)$ positively spans
  $\R^{n-d}$.\\
  {\bf (iii)} For $n\geq d\geq 2r+1$, every subset of $n-1-r$ rows of
  $\overline A(\Sigma)$ contains an alternating subset of size
  $n-d+1$.}

To see {\bf(i)}, we note that
$\{1,t,{\binom{t}{2}},\ldots,{\binom{t}{n-d+1}}\}$ and
$\{1,t,t^2,\ldots, t^{n-d+1}\}$ are two different bases for the vector
space of rational polynomials of degree at most $n-d+1$. Thus the
matrix $\overline A=\overline A(\Sigma)$ arises by invertible row
operations from the matrix
$$\overline B=((-1)^k(k-1)^{j-1})_{2\le k\le n,\ 1\le j\le n},$$
whose maximal minors are Vandermonde determinants.
\eject

For {\bf(ii)}, observe the following: It can be seen from the oriented
matroid of the vector configuration whether there is a positive linear
relation between a set of vectors, such as the rows of $\overline A$.
The condition (in order to obtain that it positively spans) is that
the configuration must be a positive circuit, or equivalently, 
{\em totally cyclic\/} \cite[Sect.~3.4]{BLSWZ}
\cite[Sect.~6.4]{Z}.  For $\varepsilon$ small enough (as indicated),
we have the same oriented matroid as for $\varepsilon=0$, and hence as
for $\overline B$. The oriented matroid ${\cal M}(\overline B)$ of
rank $n-d$ determined by the rows of $\overline B$ is the dual of the
alternating oriented matroid, which differs from the alternating
oriented matroid by a reorientation of $\{3,5,7,\ldots\}$, cf.\ 
Section~\ref{sec:om}. In ${\cal M}(\overline B)$, the positive
circuits are all subsets of size $n-d+1$. In particular, the ground
set is a positive vector, that is, there is a positive linear relation
among all the rows of~$\overline B$.

For part {\bf(iii)}, start with the alternating index set
$\{2,3,\ldots,n\}$ for the rows of $\overline A$. Now successively
delete any $r$ rows from $\overline A$, but whenever a row is deleted,
we remove also the next row above or below that has not yet been
deleted.  Thus in each of the (at most) $r$ deletion steps, we remove
two adjacent rows of $\overline A$, and hence the index set is kept to
be alternating. After all this, we are left with a submatrix
$\overline{\overline A}$ of $A$ that has at least $n-1-2r$ rows and
whose index set is alternating. Since $n-d+1\geq n-1-2r$, we may take
the first $n-d+1$ rows of $\overline{\overline A}$.  \smallskip

{\bf (C) } 
{\em The polytopes $C_d^n(\varepsilon)$ are indeed cubical.}

Let $F$ be a facet of~$C_d^n(\varepsilon)$.  Its preimage
$F'=\pi^{-1}(F)$ is a face of~$C^n(\varepsilon)$ of dimension at
least~$d-1$.  Suppose $\dim F'\ge d$.  From
Lemma~\ref{lem:orthogonal} we infer that some $n-d$ rows of the
$n\times(n-d)$-matrix~$A$ are (positively) linearly dependent.  This
is a contradiction to~{\bf(i)} of Claim~3.  We conclude that $\dim
F'=d-1=\dim F$.
The claim now follows from Lemma~\ref{lem:isomorphic}.
\endproof

It follows from the discussion in the last paragraph that each facet
of~$C_d^n(\varepsilon)$ is an affinely isomorphic image of some
$(d-1)$-face of~$C_d^n=C_d^n(\varepsilon)$.

\begin{Proposition}\label{prp:facets}~\\
  The facets of~$C_d^n$ bijectively correspond to the
  distinct positive circuits that can be found in the oriented
  matroids~${\cal M}(A(\Sigma))$ associated to the rows of the
  $n\times(n-d)$-matrices~$A(\Sigma)$, for all possible choices
  $\Sigma\in\{\pm 1\}^n$ of signs.
\end{Proposition}

An equivalent formulation of Proposition~\ref{prp:facets} is the
following: Consider the $2n\times(n-d)$-matrix $M$, whose odd numbered
rows are the rows of~$A({+}{+}\cdots{+})$, while the even numbered
rows are the rows of~$A({-}{-}\cdots{-})$.  Then the facets of~$C_d^n$
bijectively correspond to the positive circuits of~$M$ which do not
contain the first two rows, or the second two rows \ldots that is,
that do not contain both the $(2i-1)$-st and the $2i$-th row, for
any~$i$.

\begin{Theorem}{\bf (Cubical Gale Evenness Criterion)}\label{thm:cubicalGEC}\\
The facets of $C_d^n$ are given by the subsets
\[
\alpha\subseteq\{-1,+1,-2,+2,\ldots,-n,+n\},\quad\sharp\alpha=n-d+1,\quad
\alpha\cap(-\alpha)=\emptyset
\]
of the following forms:
\begin{description}
\item[\boldmath{$p=0$}:] $\pm 1\not\in\alpha$, 
and $|\alpha|$ satisfies the usual
(simplicial) Gale evenness criterion: between any two values~$\pm
i,\pm j\in\alpha$ there is an even number of ``zeroes,'' that is, an
even number of values~$k$ such that $\pm k\not\in\alpha$ for $i<k<j$.

\item[\boldmath{$1\leq p\leq n-d+1$}:]
$\alpha=\{-1,+2,\ldots,(-1)^{p-1}(p-1),\sigma\cdot p\} \cup \alpha^{(p)}$,
where $\sigma\in\{-1,+1\}$,\\
$\sharp \alpha^{(p)}=n-d+1-p=\sharp\alpha-p$, 
and if $\alpha^{(p)}\not=\emptyset$ (that is, $p\leq n-d$), then\\
{\rm(1)~} $\min \alpha^{(p)}>p+1$ [``there {\bfseries is} a gap''],\\
{\rm(2)~} $\alpha^{(p)}$
satisfies the usual (simplicial) Gale evenness criterion,
and\\
{\rm(3)~} 
if $\sigma = (-1)^{p+1}$, then 
$\min (\alpha^{(p)})-p$ is even [``the first gap is odd''];\\
\hphantom{{\rm(3)~}}
if $\sigma = (-1)^{p~~}$, then 
$\min(\alpha^{(p)})-p$ is odd [``the first gap is even''].
\end{description}
\end{Theorem}

\proof Define $p=\min\{i\geq 0:\pm (i+1)\not\in\alpha\}$. We
use this parameter to classify the positive circuits of the
$2n\times(n-d)$ matrix

$$M=\left(
\begin{array}{ccccccc}
\pm\varepsilon &&&&&&\\
1 & \pm\varepsilon &&&&&\\
-1 & -1 & \pm\varepsilon &&&&\\
1 & 2 & 1 & \ddots &&&\\
&&&&  \pm\varepsilon &&\\
&&&& (-1)^p & \pm\varepsilon &\\
&&&&& (-1)^{p} & \pm\varepsilon
\end{array}
\right)
$$

Thus $\alpha$ has the form $\{\sigma_1\cdot 1,\sigma_2\cdot
2,\ldots,\sigma_{p-1}(p-1), \sigma_p\cdot p\}\cup\alpha^{(p)}$, where
$\sigma_i\in\{+1,-1\}$. Now $\alpha^{(p)}$ must yield a positive
circuit in the contraction obtained by deleting the first $p$~rows and
columns. But this contraction is just a dual of a cyclic oriented
matroid, with every element doubled. Thus the Gale evenness condition
is both necessary and sufficient.

Next we determine the correct sign $\sigma=\sigma_p$ such that
$\{\sigma_p\cdot p\}\cup\alpha^{(p)}$ is a positive circuit of the
respective contraction, that is, so that the corresponding rows
of~$M$ have a linear combination with positive coefficients for
which the last $n-d-(p-1)$ components vanish. We get the answer
by comparison with $\{\pm(p+1)\}\cup\alpha^{(p)}$:
\begin{itemize}
\item If the gap between $(p+1)$ and $\alpha^{(p)}$ is {\em even},
then $\{\pm (p+1)\}\cup\alpha^{(p)}$ is a positive circuit of its
contraction by the Gale evenness criterion, for any small enough
$\varepsilon> 0$, and hence also for $\varepsilon=0$. Then the
``$(-1)^{p+1}$'' component of $\pm(p+1)$ can be replaced by the
``$(-1)^{p+1}\varepsilon$'' component of $(-1)^{p+1}\cdot p$. Hence\\
{\em If the gap between $p$ and $\alpha^{(p)}$ is {\bfseries odd}, then we
need $\sigma_p=(-1)^{p+1}$.}
\item If the gap between $(p+1)$ and $\alpha^{(p)}$ is {\em odd}, then
$\{\pm (p+1)\}\cup \alpha^{(p)}$ is a circuit with negative element
``$\pm(p+1)$'' and all other elements positive, since it {\em
violates} Gale's evenness criterion. But then if we replace the
``$(-1)^{p+1}$'' component of $\pm(p+1)$ by a
``$(-1)^p\varepsilon$'' component of $(-1)^p\varepsilon$, we get a
positive circuit supported on $\{(-1)^p p\}\cup\alpha^{(p)}$. Hence\\
{\em If the gap between $p$ and $\alpha^{(p)}$ is {\bfseries even}, then we
need $\sigma_p=(-1)^p$.}
\end{itemize}

The rest is easy: looking at the circuit on 
$$\{\sigma_1\cdot 1,\sigma_2\cdot 2, \ldots,\sigma_p p\}\cup
\alpha^{(p)}$$ 
that is obtained by extending the linear combination of
$\alpha^{(p)}$, we see that the coefficient of the row $\sigma_p\cdot p$ is very
large (order of $O(\frac{1}{\varepsilon})$) compared to the
coefficients on $\alpha^{(p)}$. To compensate this, we need\\
$\sigma_{p-1}= -(-1)^p$ with coefficient of order
$O(\frac{1}{\varepsilon^2})$,\\
$\sigma_{p-2}= -(-1)^{p-1}$ with coefficient of order
$O(\frac{1}{\varepsilon^3})$, and so on until\\
$\sigma_{1}= -(-1)^{2}$ with coefficient of order
$O(\frac{1}{\varepsilon^p})$.\\
The case $p=0$ is included in this.  The case
$p=n-d+1$ ($\alpha^{(p)}=\emptyset$)
is easy.
\endproof

The ``cubical Gale evenness criterion'' also
allows us to count the facets of~$C_d^n$.

\begin{Corollary}\label{cor:fnd}~~\\
  For $n\ge d\ge2$ and for $\varepsilon$ sufficiently small, then the
  number $f(n,d):=f_{d-1}(C_d^n)$ of facets
  of~$C_d^n$ is given by
\begin{eqnarray*}
f(n,d) &=& 2d \ +\  
4\sum_{p=0}^{n-d-1}\big( \binom{\lfloor\frac d2\rfloor+p+1}{p+2}+
                         \binom{\lfloor\frac{d+1}2\rfloor+p}{p+2}\big)2^p
\\
       &=& 2d \ +\  
\sum_{p=2}^{n-d+1}\big( \binom{\lfloor\frac d2\rfloor+p-1}{p}+
                         \binom{\lfloor\frac{d-1}2\rfloor+p-1}{p}\big)2^p
\\
       &=& \left\{
\begin{array}{ll}
\displaystyle{
2d \ +\  
\sum_{p=2}^{n-d+1}\binom{k+p-1}{p}2^{p+1}
}& \mbox{for $d=2k$,}\\
\displaystyle{
2d \ +\  
\sum_{p=2}^{n-d+1} \frac{p+2k-2}{p+k-1}\binom{k+p-1}{p}2^p
}& \mbox{for $d=2k+1$.}
\end{array}\right.
\end{eqnarray*}
\end{Corollary}

We evaluate the function $f(n,d)$ for some particularly interesting
choices of~$(n,d)$:
\begin{eqnarray*}
d=2: && f(n,2)\ \ =\ \  2^n\\[1mm]
d=3: && f(n,3)\ \ =\ \  2^n-2\\[1mm]
d=4: && f(n,4)\ \ =\ \  (n-2)2^{n-2}\\[1mm]
d=5: && f(n,5)\ \ =\ \  (n-4)2^{n-2}+2\\[2mm]
n=d: && f(d,d)\ \  =\ \  2d\\[1mm]
n=d+1: && f(d+1,d)\ \  =\ \  d^2+d+2\lfloor\textstyle{\frac d2}\rfloor
\end{eqnarray*}

For $d$ even, the number $f(n,d)$ is already determined by $n$
and~$d$, together with
the fact that we have a cubical $d$-polytope with the $r$-skeleton of
the $n$-cube --- using the cubical Dehn-Sommerville equations, cf.\
Proposition~\ref{prop:dehn}.

%

\section{A counter-example to the Cubical Upper Bound
  Conjecture.}\label{subsec:c46}


The following construction of a cubical $3$-sphere starts with 
the cubical $4$-polytope~$C_4^6$, whose $f$-vector is
$(64,192,192,64)$.  By local ``surgery'' (the cubical equivalent
of a bistellar flip) we will obtain a cubical $3$-sphere
with $f$-vector $(64,196,198,66)$: Thus we have a counter-example
to Kalai's cubical upper bound conjecture, verifying the second
half of Theorem~\ref{thm:n=d+1}.
The verification that the cubical flip can indeed be performed
relies heavily on the description of $C_4^6$ given by the
``cubical Gale evenness criterion''
(Theorem~\ref{thm:cubicalGEC}).
\smallskip

Adopting standard oriented matroid notation we
denote a $k$-face of the $6$-cube by a sign vector
in~$\{{+},{-},{0}\}^6$ with $k$~zeroes.  Each $k$-face of~$C_4^6$
corresponds to a $k$-face of the $6$-cube.  Thus the non-empty faces
of~$C_4^6$ also correspond to {\em certain\/} sign vectors
in~$\{{+},{-},{0}\}^6$.

We show that $C_4^6$ does not maximize the $f$-vector among cubical
$3$-spheres with $64$ vertices, as follows.  Consider the following
chain of three cubical facets $A=({-}{+}{0}{0}{+}{0})$,
$B=({-}{0}{0}{+}{+}{0})$, $C=({-}{-}{0}{+}{0}{0})$.

\begin{figure}[htbp]
  \begin{center} 
    \input{threecubes.pstex_t}
    \caption{Modifying three facets of $C_4^6$.}
{(a) The three adjacent facets that are cut out. \qquad
 (b) The ball to be glued in.}
    \label{fig:3facets}
  \end{center}
\end{figure}

Consider the cubical $3$-ball $\Phi=A\cup B\cup C$.  Its
boundary~$\partial\Phi$ is a cubical $2$-sphere.

\begin{Lemma}~~\\
  Each facet of $C_4^6$, except for $A$, $B$, or~$C$, intersects
  $\partial\Phi$ in precisely one (possibly empty) face.
\end{Lemma}

\proof
Consider the following three pairs of $2$-faces: $(A\sm B,B\sm A)$,
$(B\sm C,C\sm B)$, and $(A\sm B,C\sm B)$.  We prove that for each such
pair~$(X,Y)$ none of the vertices of~$X$ is on a common facet with any 
vertex of~$Y$.  The claim then follows. We proceed case by case.

Assume that $F$ is a facet that contains any vertex from $A\sm
B=({-}{+}{0}{-}{+}{0})$ and any vertex from $B\sm
A=({-}{-}{0}{+}{+}{0})$.  Then $F=({-}{0}{u}{0}{+}{v})$, where
$u,v\in\{{+},{-},0\}$ are to be determined.  Either $u=0$ or $v=0$.
The ``cubical Gale evenness criterion,''
Theorem~\ref{thm:cubicalGEC}, case $p=1$, implies that the initial
minus sign must be followed be an even number of zeroes, which is
impossible.

In the second case $B\sm C=({-}{+}{0}{+}{+}{0})$ and 
$C\sm B=({-}{-}{0}{+}{-}{0})$.  Hence a presumptive facet would have
coordinates $({-}{0}{w}{+}{0}{x})$ with either $w=0$ or $x=0$.  By the
``cubical Gale evenness criterion,'' Theorem~\ref{thm:cubicalGEC},
case $p=1$, the initial minus sign is followed by an even number of
zeroes, so $w=0$ and $x\in\{{+},{-}\}$.  Neither choice extends the
single zero at position~$5$ to an even number of zeroes.

In the final case we would have a facet $({-}{0}{y}{0}{0}{z})$ with
$y,z\in\{{+},{-}\}$.  But, we cannot get rid of the single zero in
position~$2$.  \endproof

Due to the preceding Lemma it is possible to replace the
subcomplex~$\Phi$ by an arbitrary cubical $3$-ball with the same
boundary without changing the topology.  In particular, the Lemma
implies that replacing $\Phi$ by a cubical $3$-ball with the same
boundary still yields a cubical complex; that is, the intersection of
any two of its faces is again a face.  One crucial (but not
sufficient) condition is that each of the eight edges in $B\sm C$ and
$B\sm A$ is contained in at least four facets.  Actually, half of them 
are contained in five facets each.

We modify the boundary complex of~$C_4^6$ by a ``local
surgery,'' as follows: remove the three facets $A$, $B$, $C$
together with the two $2$-faces between them, see
Figure~\ref{fig:3facets}(a).  Into the resulting ``hole,'' glue a cubical
ball that consists of four edges (connecting the vertices of
$A\sm B$ with the corresponding vertices of $C\sm B$), 
eight $2$-faces and five cubes (four cubes grouped around a
central $3$-cube whose top facet is in $A\sm B$, and whose
bottom facet is in~$C\sm B$).

The resulting cubical sphere $\Psi$ has the $f$-vector
\[
f(\Psi)\ \ =\ \ f(C_4^6)+(0,4,8-2,5-3).
\]
Thus $C_4^6$ is a neighborly cubical polytope whose $f$-vector is
not maximal among the cubical $3$-spheres with $64$~vertices.  We
conjecture that $\Psi$ is not polytopal.

\section{Comments.}\label{sec:comments}

{\bf (1) } A cubical $d$-polytope is {\em $k$-stacked\/} if it has a
cubical subdivision without interior $(d-k-1)$-faces
\cite[Def.~5.3]{BBC}. Thus every $k$-stacked cubical polytope is also
$(k+1)$-stacked.

\begin{Proposition}~~\\
  The neighborly cubical $d$-polytopes $C_d^n$ are
  $\lfloor\frac{d+1}{2}\rfloor$-stacked.
\end{Proposition}

\proof The projection $\R^{d+1}\rightarrow \R^d$ that deletes the last
coordinate maps $C^n_{d+1}$ to $C_d^n$. The ``upper faces'' of the
cubical polytope $C^n_{d+1}$ (those with a normal vector whose last
coordinate is positive) thus define a cubical subdivision $B^n_d$ of
$C^d_n$. Since, for $r=\lfloor\frac{d}{2}\rfloor -1$, all $r$-faces of
$C^n_{d+1}$ get mapped to the boundary of $C_d^n$, we conclude that
$C_d^n$ has no interior $r$-faces. Thus $C_d^n$ is $k$-stacked, for
$k=d-1-r=d-\lfloor \frac{d}{2}\rfloor=\lfloor \frac{d+1}{2}\rfloor$.
\endproof

In the setting of Babson, Billera \& Chan \cite[Sect.~5]{BBC} this
yields a new extreme ray for the cubical $g$-cone for even $d$.

\begin{Corollary}~~\\
  For even $d$ and $k:=\frac{d}{2}$, the neighborly cubical
  $d$-polytopes $C_d^n$ form a sequence of polytopes for which the
  cubical $g$-vector is dominated by its $k$-th component, that is
\[
\lim\limits_{n\rightarrow\infty}\frac{g^c_i(C_d^n)}{g^c_k(C_d^n)}=
0\qquad\textrm{ for all }i\not= k.
\]
\end{Corollary}

Thus ``the ray $\R^+e_k$ lies in the closure of the Adin $g$-cone,''
for even $d$ and $k=\frac{d}{2}$, in the terminology of
\cite[Sect.~5]{BBC}. In particular, for $d=4$ this implies that the
closure of the Adin $g$-cone, i.e., of
\[
\textrm{cone\,} \{(g^c_1, g^c_2)\ =\ 
(f_0-16, 4f_3-3f_0+16):\textrm{ cubical $4$-polytopes}\} 
\]
is the complete positive orthant in~$\R^2$.
\medskip

{\bf (2) } One is tempted to ask whether for {\em any\/}
polytope $P$ of dimension $n\ge d$ there is a $d$-polytope $Q$ that
has isomorphic $r$-skeleton, for $r=\lfloor d/2\rfloor -1$.  In
general, this is known to be false, due to a construction by Klee and
Gr\"unbaum, see Gr\"unbaum~\cite[12.2.2]{Gr}.  But the case of the
cross polytopes, $P=C_n^{\triangle}$, seems to be open and
particularly interesting.  \medskip
\medskip

{\bf (3) } Is there a construction of (even-dimensional) neighborly
cubical polytopes that have all vertices on a sphere?  Note that the
trigonometric moment curve \cite[p.~67]{Gr} \cite[p.~75]{Z} yields
this in the simplicial case.  By explicit construction of Schlegel
diagram as the Delaunay subdivisions of a finite point set, Raimund
Seidel~\cite{S} obtained such a construction for $d=4$ and $n\le7$.
\medskip

{\bf (4) } The polytope $C_4^5\cong P$ and the non-cubical polytope~$P'$ of
Section~\ref{subsec:notcubical},
as well as the polytope $C_4^6$ and 
the counter-example of Section~\ref{subsec:c46}
were first constructed using the {\sc polymake} system~\cite{GJ}. 
Corresponding {\sc polymake} input files are available 
at

\begin{quote} {\tt
  http://www.math.tu-berlin.de/diskregeom/polymake/examples/NCP/}.
\end{quote}

The sequence of facets $A$, $B$, $C$ used in Section~\ref{subsec:c46}
was also found in an exhaustive computer
search.  The same test revealed that a similar construction does not
work with the smaller $C_4^5$ polytope.
\bigskip

\vfill
\noindent{\bf Acknowledgement.}
Thanks to Raimund Seidel for inspiring e-discussions.  We are\break indebted
to Roswitha Blind for a correction that led to a substantial
modification in the construction of
the counter-example of Section~\ref{subsec:c46}.
\eject

\end{document}

%% file: quarter.pstex_t
\begin{picture}(0,0)%
\epsfig{file=quarter.pstex}%
\end{picture}%
\setlength{\unitlength}{1973sp}%
\begingroup\makeatletter\ifx\SetFigFont\undefined%
\gdef\SetFigFont#1#2#3#4#5{%
  \reset@font\fontsize{#1}{#2pt}%
  \fontfamily{#3}\fontseries{#4}\fontshape{#5}%
  \selectfont}%
\fi\endgroup%
\begin{picture}(4785,5680)(1377,-6661)
\put(1701,-1161){\makebox(0,0)[rb]{\smash{\SetFigFont{7}{8.4}{\rmdefault}{\mddefault}{\updefault}$5$}}}
\put(1701,-2161){\makebox(0,0)[rb]{\smash{\SetFigFont{7}{8.4}{\rmdefault}{\mddefault}{\updefault}$4$}}}
\put(1701,-3161){\makebox(0,0)[rb]{\smash{\SetFigFont{7}{8.4}{\rmdefault}{\mddefault}{\updefault}$3$}}}
\put(2776,-1861){\makebox(0,0)[rb]{\smash{\SetFigFont{7}{8.4}{\rmdefault}{\mddefault}{\updefault}$\alpha$}}}
\put(1701,-4161){\makebox(0,0)[rb]{\smash{\SetFigFont{7}{8.4}{\rmdefault}{\mddefault}{\updefault}$2$}}}
\put(1701,-5161){\makebox(0,0)[rb]{\smash{\SetFigFont{7}{8.4}{\rmdefault}{\mddefault}{\updefault}$1$}}}
\put(1701,-6161){\makebox(0,0)[rb]{\smash{\SetFigFont{7}{8.4}{\rmdefault}{\mddefault}{\updefault}$0$}}}
\put(2001,-6661){\makebox(0,0)[b]{\smash{\SetFigFont{7}{8.4}{\rmdefault}{\mddefault}{\updefault}$0$}}}
\put(3001,-6661){\makebox(0,0)[b]{\smash{\SetFigFont{7}{8.4}{\rmdefault}{\mddefault}{\updefault}$1$}}}
\put(4001,-6661){\makebox(0,0)[b]{\smash{\SetFigFont{7}{8.4}{\rmdefault}{\mddefault}{\updefault}$2$}}}
\put(5001,-6661){\makebox(0,0)[b]{\smash{\SetFigFont{7}{8.4}{\rmdefault}{\mddefault}{\updefault}$3$}}}
\put(6001,-6661){\makebox(0,0)[b]{\smash{\SetFigFont{7}{8.4}{\rmdefault}{\mddefault}{\updefault}$4$}}}
\end{picture}

%% file: threecubes.pstex_t
\begin{picture}(0,0)%
\epsfig{file=threecubes.pstex}%
\end{picture}%
\setlength{\unitlength}{1579sp}%
\begingroup\makeatletter\ifx\SetFigFont\undefined%
\gdef\SetFigFont#1#2#3#4#5{%
  \reset@font\fontsize{#1}{#2pt}%
  \fontfamily{#3}\fontseries{#4}\fontshape{#5}%
  \selectfont}%
\fi\endgroup%
\begin{picture}(18910,9240)(-21,-8386)
\put(4351,-5011){\makebox(0,0)[rb]{\smash{\SetFigFont{10}{12.0}{\familydefault}{\mddefault}{\updefault}3}}}
\put(5626,-4561){\makebox(0,0)[rb]{\smash{\SetFigFont{10}{12.0}{\familydefault}{\mddefault}{\updefault}$({-}{-}{-}{+}{+}{+})$}}}
\put(5101,-5536){\makebox(0,0)[b]{\smash{\SetFigFont{10}{12.0}{\familydefault}{\mddefault}{\updefault}6}}}
\put(4426,-6511){\makebox(0,0)[rb]{\smash{\SetFigFont{10}{12.0}{\familydefault}{\mddefault}{\updefault}5}}}
\put(4126,-4036){\makebox(0,0)[lb]{\smash{\SetFigFont{10}{12.0}{\familydefault}{\mddefault}{\updefault}2}}}
\put(4501,-1486){\makebox(0,0)[rb]{\smash{\SetFigFont{10}{12.0}{\familydefault}{\mddefault}{\updefault}4}}}
\put(9825,-2461){\makebox(0,0)[lb]{\smash{\SetFigFont{10}{12.0}{\familydefault}{\mddefault}{\updefault}$({-}{+}{-}{+}{+}{-})$}}}
\put(6750,539){\makebox(0,0)[b]{\smash{\SetFigFont{10}{12.0}{\familydefault}{\mddefault}{\updefault}$A\sm B=({-}{+}{0}{-}{+}{0})$}}}
\put(7425,-661){\makebox(0,0)[lb]{\smash{\SetFigFont{10}{12.0}{\familydefault}{\mddefault}{\updefault}$({-}{+}{+}{-}{+}{-})$}}}
\put(8625,-61){\makebox(0,0)[lb]{\smash{\SetFigFont{10}{12.0}{\familydefault}{\mddefault}{\updefault}$({-}{+}{-}{-}{+}{-})$}}}
\put(5100,-7861){\makebox(0,0)[rb]{\smash{\SetFigFont{10}{12.0}{\familydefault}{\mddefault}{\updefault}$({-}{-}{+}{+}{-}{+})$}}}
\put(8625,-7411){\makebox(0,0)[lb]{\smash{\SetFigFont{10}{12.0}{\familydefault}{\mddefault}{\updefault}$({-}{-}{-}{+}{-}{-})$}}}
\put(7350,-7861){\makebox(0,0)[lb]{\smash{\SetFigFont{10}{12.0}{\familydefault}{\mddefault}{\updefault}$({-}{-}{+}{+}{-}{-})$}}}
\put(9825,-4861){\makebox(0,0)[lb]{\smash{\SetFigFont{10}{12.0}{\familydefault}{\mddefault}{\updefault}$({-}{-}{-}{+}{+}{-})$}}}
\put(7950,-3136){\makebox(0,0)[lb]{\smash{\SetFigFont{10}{12.0}{\familydefault}{\mddefault}{\updefault}$({-}{+}{+}{+}{+}{-})$}}}
\put(7950,-5536){\makebox(0,0)[lb]{\smash{\SetFigFont{10}{12.0}{\familydefault}{\mddefault}{\updefault}$({-}{-}{+}{+}{+}{-})$}}}
\put(5100,-436){\makebox(0,0)[rb]{\smash{\SetFigFont{10}{12.0}{\familydefault}{\mddefault}{\updefault}$({-}{+}{+}{-}{+}{+})$}}}
\put(6300, 14){\makebox(0,0)[rb]{\smash{\SetFigFont{10}{12.0}{\familydefault}{\mddefault}{\updefault}$({-}{+}{-}{-}{+}{+})$}}}
\put(6825,-8386){\makebox(0,0)[b]{\smash{\SetFigFont{10}{12.0}{\familydefault}{\mddefault}{\updefault}$C\sm B=({-}{-}{0}{+}{-}{0})$}}}
\put(5625,-2236){\makebox(0,0)[rb]{\smash{\SetFigFont{10}{12.0}{\familydefault}{\mddefault}{\updefault}$({-}{+}{-}{+}{+}{+})$}}}
\put(10050,-1411){\makebox(0,0)[lb]{\smash{\SetFigFont{10}{12.0}{\familydefault}{\mddefault}{\updefault}$A=({-}{+}{0}{0}{+}{0})$}}}
\put(10050,-6286){\makebox(0,0)[lb]{\smash{\SetFigFont{10}{12.0}{\familydefault}{\mddefault}{\updefault}$C=({-}{-}{0}{+}{0}{0})$}}}
\put(6225,-7111){\makebox(0,0)[rb]{\smash{\SetFigFont{10}{12.0}{\familydefault}{\mddefault}{\updefault}$({-}{-}{-}{+}{-}{+})$}}}
\put(3825,-5236){\makebox(0,0)[rb]{\smash{\SetFigFont{10}{12.0}{\familydefault}{\mddefault}{\updefault}$({-}{-}{+}{+}{+}{+})$}}}
\put(3825,-2911){\makebox(0,0)[rb]{\smash{\SetFigFont{10}{12.0}{\familydefault}{\mddefault}{\updefault}$({-}{+}{+}{+}{+}{+})$}}}
\put(10125,-3811){\makebox(0,0)[lb]{\smash{\SetFigFont{10}{12.0}{\familydefault}{\mddefault}{\updefault}$B=({-}{0}{0}{+}{+}{0})$}}}
\end{picture}